\documentclass[a4paper,12pt]{article}
\usepackage{amssymb}
\usepackage{latexsym}
\usepackage{amsmath}
\usepackage{amsthm}
\usepackage{amsfonts}
\usepackage{amssymb}
\usepackage{graphicx}

\newtheorem{Th}{Theorem}

\newcommand{\be}{\begin{equation}}
\newcommand{\ee}{\end{equation}}
\newcommand{\R}{\mathbb{R}}
\newcommand{\N}{\mathbb{N}}
\newcommand{\C}{\mathbb{C}}

\newcommand{\reset}{\setcounter{equation}{0}\setcounter{Th}{0}\setcounter{Prop}{0}\setcounter{Co}{0}
\setcounter{Lm}{0}\setcounter{Rm}{0}}
\def\La{\Lambda}
\def\La{\Lambda}

\def\lf{\left}
\def\rg{\right}

\def\al{\alpha}
\def\la{\lambda}

\def\ds{\displaystyle}
\def\ov{\overline}
\def\Om{\Omega}

\def\p{\partial}

\begin{document}
\title{Exploring the Unknown: the Work of Louis Nirenberg on Partial Differential Equations.}
\author{ Tristan Rivi\`ere\footnote{Department of Mathematics, ETH Z\"urich,
CH-8093 Z\"urich, Switzerland.}}
\maketitle
\section{Preamble}
\reset
Partial differential equations are central objects in the mathematical modeling of natural and social sciences (sound propagation, heat diffusion, thermodynamics, electromagnetism, elasticity, fluid dynamics, quantum mechanics, population growth, finance...etc). They were for a long time restricted only to the study of natural phenomena or questions pertaining to geometry, before becoming over the course of time, and especially in the last century, a field in itself. 

The second half of the XXth century was the ``golden age'' of the exploration of partial differential equations from a theoretical perspective. 

The mathematical work of Louis Nirenberg since the early 1950s has to a large extent contributed to the growth of this fundamental area of human knowledge. The name Nirenberg is associated with many of the milestones in the study of PDEs. 

The award of the {\it Abel Prize} to Louis Nirenberg marks a special occasion for us to revisit the development of the field of PDEs and the work of one of the main actors of its exploration.

\section{Introduction}
\reset

Louis Nirenberg uses to qualify the field of {\it Partial Differential Equations} as being ''messy'' (and also often acknowledges his special taste for what he called  this ''messiness'' )
which is probably a reference to the ''intrinsic diversity'' of the field.
We would like to illustrate the pertinence of  this quote by contradiction and by presenting the original attempts made mostly in the XIXth century
to see PDE as a whole and the limits and inadequacies this approach has been confronted to.

\subsection{A general existence result: the Cauchy-Kowalevski theorem}
Perhaps the first general systematic study of partial differential equations goes back to the work of Augustin-Louis Cauchy and his existence theorem for
{\it quasilinear} first order PDEs with real analytic data.
\begin{Th}
\label{th-I.1} {\bf[Cauchy 1842]} Let $(A^j)_{j=1\cdots {n-1}}$ be a family of $n-1$ real analytic maps from a neighborhood of $(0,0)$ in ${\R}^n\times{\R}^m$ into the space of real $m\times m$ matrices
and let $f$ be a real analytic map into ${\R}^m$. Then there exists a unique real analytic solution $u$ in a neighborhood of the origin to the following system
\[
\lf\{
\begin{array}{l}
\ds\p_{x_n}u=\sum_{j=1}^{n-1} A^j(x,u)\ \p_{x_j}u+f(x,u)\\[5mm]
\ds u(x_1,\cdots,x_{n-1},0)\equiv 0\quad.
\end{array}
\rg.
\]
\hfill $\Box$
\end{Th}
In 1874, Sofia Kowalevski (or Kowalevskaya), apparently unaware of Cauchy's work, proved in her thesis a general nonlinear version of the previous result. We present the theorem called nowadays {\it Cauchy-Kowalevski theorem} in the particular case of {\it second order non-linear} scalar equations
\[
F(x,u,\p u,\p^2u)=0\quad\mbox{where }\quad\p u:=(\p_{x_i}u)_{i=1\cdots n}\quad\mbox{ and }\quad\p^2 u:=(\p^2_{x_i x_j} u)_{i,j=1\cdots n}\quad.
\]
Here $F(x,q,p,r)$ is a real analytic function of all the entries $p\in {\R}$, $q=(q_i)_{i=1\cdots n}\in{\R}^n$ and $r=(r_{ij})_{i,j=1\cdots n}\in{\R}^{n^2}$. To that purpose, we introduce the notion of {\it characteristic direction}. A direction $X=(\xi_i)_{i=1\cdots n}\in {\R}^n$ is called {\it characteristic} at $(x,p,q,r)$ if
\[
\sum_{i,j}^n\frac{\p F}{\p r_{ij}}(x,p,q,r)\,\xi_i\,\xi_j=0\quad.
\]
Given a function $u$ defined in a neighborhood of a point $x_0$, an hyper-surface $f(x)=0$ is called {\it non-characteristic at $x_0$} if $\xi:=\nabla f(x_0)$ is not a characteristic direction at $(x_0,u(x_0),\p u(x_0),\p^2u(x_0))$. The Cauchy-Kowaleski theorem requires that {\it initial data}, the so called {\it Cauchy data}, be given on a {\it non-characteristic} hyper-surface. Observe that if the PDE has the form
\[
F(x,u,\p u,\p^2u)=\p^2_{x_n x_n}u-f(x,u,\p u, \p^2 u)\:,
\]  
where $f$ does not depend on $r_{nn}$, the surface $x_n=0$ is automatically {\it non characteristic}. We can now state the {\it Cauchy-Kowalevski theorem} in that particular case.
\begin{Th}
\label{th-I.2} {\bf [Cauchy-Kowalevski 1874]} Let $f(x,p,q,r)$ be a real analytic function of its variables, and assume that $f$ is independent $r_{nn}$. Let $u_0(x_1,\cdots,x_{n-1})$ and 
$u_1(x_1,\cdots,x_{n-1})$ be two real analytic functions defined in the neighborhood of the origin $0$. Then there is a unique real analytic solution defined in a neighborhood of the origin for the problem
\[
\lf\{
\begin{array}{l}
\ds \p^2_{x_n x_n}u=f(x,u,\p u, \p^2 u)\\[5mm]
\ds u=u_0\quad\mbox{ and }\quad \p_{x_n}u=u_1\quad\mbox{ on }\quad x_n=0\quad.
\end{array}
\rg.
\]
\hfill $\Box$
\end{Th}
\subsection{ Some inadequacies of the Cauchy-Kowalevski theory}
The Cauchy-Kowalevski theorem requires an analytic framework. Its proof (the historical one), consists of an argument based on the convergence of power series. Only the analyticity assumptions with respect to $x_n$ could be relaxed. The question of whether there could be more solutions for the same analytic data (in the $C^\infty$ class for instance) has stimulated much research, and, although there are uniqueness theorems for some classes of linear PDEs, there are also counter-examples to uniqueness (see \cite{Met}). The general question remains to be settled.  

\medskip
If one seeks global solutions, which are expected to exist in physical problems, there is an ``intrinsic'' need to relax the analytic framework, since singularities can appear in ``finite time'' (the time variable here being $x_n$) even though all data are analytic. Consider for instance the Cauchy solution to 
\[
\lf\{
\begin{array}{l}
\ds \p_{x_2}u=u\,\p_{x_1}u\\[5mm]
\ds u(x_1,0)=-\frac{x_1}{1+x_1^2}\quad.
\end{array}
\rg.
\]
The gradient on each level set of $u=u_0=-\frac{x_0}{1+x_0^2}$ for any $u_0\in(-1,1)$ is constant, hence each level set is made of straight segments leaving the point $(x_0,0)$ in the direction given by the vector 
 $(1,u_0)$. These segments have to meet at points where $u$ thus necessarily ceases to be continuous.\\

In her thesis, Sofia Kowalevski illustrated the need for the assumption that the initial surface be {\it non-characteristic} by the following example. Consider locally
an analytic solution of
\[
\lf\{
\begin{array}{l}
\ds \p_{x_2}u=\p^2_{x_1 x_1}u\\[5mm]
\ds u(x_1,0)=\frac{1}{1-x_1}\quad.
\end{array}
\rg.
\]
An elementary computation gives the explicit value of the coefficients of the Taylor series expansion for the solution, which happens not to converge at the origin. Hence, near the origin, there is no analytic solution to that equation with the given initial data. The reason why Cauchy's theorem does not apply in this case, is that $x_2=0$ is {\it characteristic} for the equation.
However this equation is nothing but the {\it heat equation} modeling in particular the diffusion of heat in an homogeneous material starting from the data $u(x_1,0)$. A non-existence result certainly seems counter-intuitive! It suggests that one should leave the analytic framework imposed by Cauchy and Kowalevski.\\

{The conclusion at this stage is that the analysis of PDEs \underbar{cannot be captured in a single theory} having simply to do with the convergence of power series.}

\section{Local solvability}
\reset

\subsection{The notion of local solvability and Lewy's counterexample}
The first attempt to go beyond the Cauchy-Kowalevski theory would be to give up uniqueness requirements and look at ``germs'' of PDEs at a point. We are thus led to considering the linear framework and we ask whether one can enlarge the class of possible solutions from the analytic class to the $C^\infty$ class, or even to the space of distributions\footnote{We return to this notion of weak solutions introduced by Laurent Schwartz later in this presentation.} ${\mathcal D}'$. Is a linear PDE always {\it locally} solvable? Precisely, consider $a_\al(x)$ to be a family of $C^\infty$ complex coefficients defined in a neighborhood of the origin in ${\R}^n$ and indexed by $\al=(\al_1,\cdots\al_n)$ where $\al_j\in {\N}$  and $|\al|:=\sum\al_j\le p$, for some $p\in {\N}^\ast$. Given any $f \in C^\infty_0({\R}^n,{\C})$, does there exist locally near the origin a complex-valued solution (if need be, even in the very weak sense of distributions) to
\[
L\,u=\sum_{\al=(\al_1\cdots\al_n)} a_{\al}(x)\,\frac{\p^{|\al|}u}{\p ^{\al_1}x_1\cdots\p ^{\al_n}x_n}=f\quad\quad ?
\]
If for a given $L$ and any $f$, the answer to the question is ``yes'', the PDE is called {\it locally solvable} at 0. Around 1955, Leon Ehrenpreis and Bernard Malgrange proved independently
the {\it local solvability} of any linear PDE with constant coefficients $a_\al(x)\equiv a_\al$ (see \cite{Ehr}, \cite{Mal}). Using Laurent Schwartz theory of tempered Distribution such a PDE can be converted into a convolution equation, and the problem is reduced to a division problem in function algebra. Encouraged by this result, the conjecture asserting that any linear operator $L$ should be {\it locally solvable} became notorious in the PDE community. That was until Hans Lewy came up 
in 1957 with a spectacular counter-example, namely: there exists a non-analytic function $f\in C^\infty$ such that the PDE
\[
\frac{\p u}{\p {x_1}}+i\,\frac{\p u}{\p {x_2}}-2\, i\ (x_1+i\,x_2)\ \frac{\p u}{\p x_3}=f
\]
has no $C^1$ solution in any neighborhood of the origin in ${\R}^3$ (see \cite{Lew}). This counter-example triggered an intense research activity which involved many prominent analysts of that time (such as Lars H\"ormander, Louis Nirenberg, Fran\c cois Tr\`eves, Yuri Vladimirovich Egorov,...). They were seeking what necessary and sufficient conditions would ensure {\it local solvability}. 

\subsection{The Nirenberg-Tr\`eves local solvability condition}

In three fundamental papers, Louis Nirenberg and Fran\c cois Tr\`eves (see \cite{NT1}, \cite{NT2} and \cite{NT3}) in 1963 have proposed a condition on the  {\it principal symbol of the PDE}
\[
p_m(x,\xi):=\sum_{|\al|=p}\ a_\al(x)\ \prod_{j=1}^n\xi_j^{\al_j}
\]
that should be necessary and sufficient for local solvability, the so called {\it (P) condition}. The bi-characteristics of a real function $A(x,\xi)$ are the curves $(x(t),\xi(t))$ that solve the Hamilton-Jacobi equation
\[
\frac{d x}{dt}=\nabla_\xi A(x,\xi)\quad\mbox{ and }\quad\frac{d \xi}{dt}=-\nabla_xA(x,\xi)\quad.
\]  
On such a curve, $A(x(t),\xi(t))$ is constant, and the curves for which this constant is zero are called the {\it null bi-characteristics} of $A$. The {\it Nirenberg-Tr\`eves (P) condition} reads as follows
\[
(P) \quad\mbox{\it On every null bi-characteristics of }{\Re}\, p_m\:\mbox{\it the function }{\Im}\, p_m\:\mbox{\it does not change sign.}
\]
The necessity and sufficiency of the {\it Nirenberg-Tr\`eves $(P)$ condition} for the {\it local solvability} was established for a growing number of cases in successive works by Nirenberg and Tr\`eves themselves, 
then by Beals and Fefferman \cite{BF}, by H\"ormander \cite{Ho2}, until a very recent result of Nils Dencker in 2006 (see \cite{Den}) which establishes the sufficiency of the generalized {\it Nirenberg-Tr\`eves $(P)$ condition} for  general {\it pseudo-differential operators}\footnote{This generalized notion of differential operators was introduced by Louis Nirenberg in collaboration with Joseph J. Kohn in \cite{KN}, and by Lars H\"ormander in \cite{Ho} few years after the introduction by Calder\'on and Zygmund of the theory of {\it singular integral operators}. This revival of the Calder\'on Zygmund theory  has probably be strongly stimulated by the proof of the index theorem by Michael Atiyah and Isadore Singer in 1962 .}, the so called {\it Nirenberg-Tr\`eves $(\psi)$ condition}. 
\section{Cauchy problems and global solvability for linear PDEs}
\reset
\subsection{The notion of Cauchy problem}
Understanding the local solvability of a PDE is certainly one important question, but one might argue that its relevance should not come into play in physical applications where solutions are expected to be {\it globally} defined on a given subdomains of spacetime. \\

The question of global solvability is traditionally coupled with that of uniqueness, and together they form what is called a {\it Cauchy problem}. A {\it Cauchy problem}, or {\it well-posed problem}, consists of a linear PDE $L$, of a function space $E$ in which the data (the input) makes sense, of a function space $F$ to which the expected solution (the output -- also called the ``unknown'') should belong, with the requirements that:
\begin{itemize}
\item[i)] there exists exactly one solution in $F$ of the PDE $L$ for any given data in the function space $E$ ;
\item[ii)] the dependence of the solution on the data is continuous from $E$ into $F$.
\end{itemize}
Finding appropriate {\it Cauchy problems} for linear PDEs has generated a tremendous amount of mathematical activities in the last century. It also gives the ``asymptotes'', or ``constraints at the horizon'', in order to solve many {\it non-linear} problems, as we shall see.\\ 

One can illustrate the notion of {\it Cauchy problem} by looking at examples of {\it ill-posed problems}. One of them deals with finding a holomorphic extension in the disc $D^2$ of a prescribed $C^1$ boundary data. In other words, given $\phi\in E:=C^1(\p D^2,{\C})$, one seeks $u\in F:=C^{1}(\ov{D^2},{\C})$ satisfying
\[
\lf\{
\begin{array}{l}
\ds Lu:=\ov{\p}u=\frac{1}{2}\lf(\frac{\p u}{\p x_1}+\, i\, \frac{\p u}{\p x_2}\rg)=0\quad\quad \mbox{ in }D^2\\[5mm]
u=\phi\quad\quad\mbox{ on }\p D^2\quad.
\end{array}
\rg.
\]
This is an ill-posed problem, because not every $\phi\in C^1( \p D^2,{\C})$ admits a $C^1$ holomorphic extension in $D^2$ (take for instance $\phi(\theta)=e^{-i\, \theta}$). 
The ill-posedness can however be thwarted by replacing $E:=C^1(\p D^2,{\C})$ with the subspace of $C^1$ functions that are $L^2$-orthogonal to $e^{-i\,k\,\theta}$ for any $k\in {\N}^\ast$.

\subsection{The fragmentation of the analysis of PDEs}

The search for Cauchy problems has imposed a fragmentation of the field of PDEs into multiple areas of analysis, which have often been developed independently of each other. This is due to the very different behaviors of various differential operators. The analysis of PDEs cannot be encapsulated into a single theory, and it is intrinsically split. It is in fact a field that might seem disorderly from the outside. But this ``messiness", which sometimes discourages young vocations, is in the very nature of PDEs and it is the source of an infinite diversity of phenomena, arguments, and results. One may nevertheless attempt to put some order in this diversity, by singling out three main families of operators that we briefly present in the framework of {\it linear second order scalar equations} with constant coefficients:
\[
Lu:=\sum_{i,j=1}^n a_{ij}\, \frac{\p^2 u}{\p x_i\, \p\, x_j}+\sum_{i=1}^n\, b_i\, \frac{\p u}{\p x_i}+ c\,u
\]
with symmetric principal symbol : $a_{ij}=a_{ji}$.
\begin{itemize}

\item[i)] The {\it elliptic operators} are those for which $L$ has no characteristic direction:
\[
\forall \xi\ne 0\quad\quad \sum_{i,j=1}^{n} a_{ij}\,\xi_i\,\xi_j>0\quad.
\]
The archetype of the elliptic operator is the {\it laplacian operator}
\[
Lu:=\Delta u=\sum_{i=1}^n \frac{\p^2 u}{\p x_i\p x_i}\quad.
\]
The {\it Laplace equation} $\Delta u=0$, also called the {\it potential equation}, has been introduced in the physics of gravitational and electromagnetic fields. The field is generated
by the potential $u$ satisfying $\Delta u=f$, where $f$ denotes either the mass or the charge density.
\item[ii)] The {\it parabolic operators} are those for which one direction is singled out and does not appear in the second derivatives. They have the form
\[
Lu:=\frac{\p u}{\p t}-\sum_{i,j=1}^{n-1} a_{ij}\, \frac{\p^2 u}{\p x_i\, \p\, x_j}-\sum_{i=1}^{n-1}\, b_i\, \frac{\p u}{\p x_i}- c\,u\quad,
\]
where the $(a_{ij})$ operator for the remaining variables is (negative) elliptic:
\[
\forall \xi\ne 0\quad\quad \sum_{i,j=1}^{n-1} a_{ij}\,\xi_i\,\xi_j>0\quad.
\]
The model case of a parabolic operator is the {\it heat operator}
\[
Lu:=\frac{\p u}{\p t}-\Delta u=\frac{\p u}{\p t}-\sum_{i=1}^{n-1} \frac{\p^2 u}{\p x_i\p x_i}\quad.
\]
The {\it heat equation} $\p_tu-\Delta u=0$ is a fundamental tool in the modeling of diffusion (heat, combustion, chemical reactions, ...).
\item[iii)] The {\it hyperbolic operators} are those for which the matrix $(a_{ij})$ has signature $(+,-,\cdots,-)$.
The
prototype is the {\it wave operator}
\[
Lu:=\frac{\p^2 u}{\p t^2}-\Delta u=\frac{\p^2 u}{\p t^2}-\sum_{i=1}^{n-1} \frac{\p^2 u}{\p x_i\p x_i}\quad.
\]
The {\it wave equation} $\p^2_{t^2}u-\Delta u=0$ is a fundamental equation arising in the modelling of propagation of waves in gases and fluids.

\end{itemize}
This (overly) simplified classification leaves out many equations, including ones that are relevant to physical phenomena, such as the Schr\"odinger equation or the Korteweg-deVries water-wave equation. Nonetheless, understanding how the three basic families differ from one another constitutes a first step in the study of ``hybrid" and of more complicated PDEs.  Parabolic equations can be understood as elliptic equations with time propagation, and thus these two families enjoy many similar properties, such as smoothing effect and infinite speed propagation. Hyperbolic equations on the other hand are very different, and they involve, for examples, finite speed propagation and the ``transport'' of singularities. The work of Louis Nirenberg has a barycenter whose location lies closer to the first two families, and this is why we shall mostly restrict this discussion to {\it elliptic} and {\it parabolic} PDEs.

\subsection{The birth of Cauchy problems: the Dirichlet principle} 

Peter Gustav Lejeune Dirichlet was a German mathematician of the first half of the XIXth century. His interest for physics led him to formulate a mathematical
problem related to potential theory, a field which finds its roots in Newtonian mechanics. This problem, named after him by Riemann, is formulated as follows:

\medskip

{\it  A bounded domain empty of charge being given as well as the value of the potential at its boundary, show mathematically that it produces a unique electromagnetic field in the interior of the domain.}

\medskip
In his lecture\footnote{published after his death by F. Grube in 1876, see \cite{Dir}.} on potential theory in G\"ottingen in 1857, Dirichlet gave a more precise mathematical formulation of the problem. Namely, show that for any bounded connected (smooth) domain $\Omega$ in ${\R}^3$ and for any continuous differentiable function $\phi$ on $\p \Omega$, there exists a unique continuously differentiable function $u$ equal to $\phi$ on $\p \Omega$ and minimizing the following integral, now called {\it Dirichlet energy},
\[
E(u):=\int_\Omega|\nabla u|^2\ dx^3\quad.
\]
Such a function satisfies the Laplace equation
\be
\label{I.1}
\lf\{
\begin{array}{l}
\ds\Delta u=0\quad\quad\mbox{ in }\Omega\\[5mm]
\ds u=\phi \quad\quad\mbox{ on }\p\Omega\quad.
\end{array}
\rg.
\ee
The claims made by Dirichlet in his lecture were not completely settled, but they were however used in the following years by Bernhard Riemann in two-dimensional complex analysis. This situation created some trouble in the mathematical community and the Dirichlet principle came under attack, in particular in the work of Karl Weierstrass, who produced examples of coercive Lagrangians (akin to the Dirichlet energy) which do not possess  minimizers for some boundary data.
The revival of the Dirichlet principle came in the early XXth century in the work of David Hilbert who solved it rigorously in some special cases and devoted the 20th of his famous list of 23 problems
to this principle. It is hard to find papers of the beginning of the XXth century on PDE without any mention of the {\it Dirichlet principle}. A detailed description of the steps towards the full resolution of Dirichlet's assertions would take us on too long, although very instructive, a detour. The resolution of Dirichlet's problem gave rise to numerous ad-hoc tools from functional analysis for solving PDEs, such as the notion of {\it Sobolev spaces}, {\it distributions}, {\it weak solutions}, ... Let us also mention the names of Henri L\'eon Lebesgue, Leonida Tonelli, Sergei Sobolev, Laurent Schwartz. \\

The modern way to solve the Dirichlet principle goes as follows. One introduces the so called {\it Sobolev Space} of Lebesgue measurable functions $u$ on the bounded domain $\Omega$ whose Schwartz distributional derivatives lie in the {\it Hilbert space} of square integrable functions:
\[
W^{1,2}(\Omega,{\R})=\lf\{u\in L^2\quad;\quad \partial_{x_i}u\in L^2\quad\mbox{for }i=1,\ldots, 3\rg\}\quad,
\]
where the distributional derivatives of $u$ are defined by duality with compactly supported smooth functions
\[
\forall\phi\in C^\infty_0(\Omega,{\R}) \quad\quad\lf<\p_{x_i}u,\phi\rg>:=-\int_{\Omega}u\, \frac{\p \phi}{\p x_i}\ dx^3\quad.
\]
This space is equipped with the scalar product for which it is complete, namely,
\[
(u,v):=\int_\Om \sum_{i=1}^3 \p_{x_i}u\,\p_{x_i}v+u\, v\ dx^3\quad.
\]
The Hilbert Sobolev space $W^{1,2}$ has favorable analytical properties, such as  {\it weak sequential pre-compactness} (from any sequence with uniformly bounded norm one can extract  
a subsequence converging weakly -- in duality with smooth functions -- to a function in $W^{1,2}$). Furthermore, the Dirichlet energy $E(u)$ is {\it weakly sequentially lower semi-continuous} and for any weakly converging sequence $u_k$  with limit $u_\infty$, there holds
\[
\liminf_{k\rightarrow +\infty} E(u_k)\ge E(u_\infty)\quad.
\]
One proves that smooth functions are dense in this {\it Hilbert space}. Moreover, the map which to a smooth function in $\ov{\Om}$ assigns its restriction to the boundary $\p\Om$ is continuous  from $W^{1,2}(\Omega,{\R})$ into the space of square integrable functions on the boundary $\p \Om$, and therefore extends continuously to a well defined map from $W^{1,2}(\Omega,{\R})$ into $L^2(\p\Om)$. These facts combined together
give the existence of a minimizer of the Dirichlet energy for any weak boundary data ($L^2$ function $\phi$) which admits a $W^{1,2}$ extension inside $\Omega$. This subspace of $L^2(\p \Om)$
is denoted by $H^{1/2}(\p\Om)$. The use of distribution theory makes it possible to assert that any such minimum $u$ satisfies the Laplace equation in a weak sense:
 \[
 \forall\phi\in C^\infty_0(\Omega,{\R}) \quad\quad\lf<\Delta u,\phi\rg>:=\int_{\Omega}u\, \Delta\phi\ dx^3=0\quad.
 \]
 A fundamental lemma of Laurent Schwartz states that such a function $u$ has to be analytic in the interior of the domain $\Omega$ and is the unique solution to (\ref{I.1}). This full resolution of the Cauchy problem for the Dirichlet question can be extended without much difficulty in the case when ``inside charges'' are present, and one shows that for any smooth function $f$ in $\Omega$ and any smooth function $\phi$ in $\p\Om$, there is a unique solution -- smooth inside $\Om$ -- to the inhomogeneous problem
 \be
\label{I.2a}
\lf\{
\begin{array}{l}
\ds-\Delta u=f\quad\quad\mbox{ in }\Omega\\[5mm]
\ds u=\phi \quad\quad\mbox{ on }\p\Omega\quad.
\end{array}
\rg.
\ee
\subsection{The Agmon-Douglis-Nirenberg elliptic Cauchy problems in the Banach $L^p$ spaces}

Numerous problems from classical mechanics, quantum mechanics, chemistry, biology, and from geometry can be described by means of elliptic non-linear PDEs of
the form
\be
\label{I.3a}
-\Delta u= f(x,u,\nabla u)\quad,
\ee
where the highest-order term is given by the Laplace operator or more generally by an {\it elliptic operator} which has no null characteristic.
Considering this PDE on a bounded smooth domain $\Omega$ with zero boundary condition, it is tempting to view it as a perturbation of (\ref{I.2})
and use the inverse of the Dirichlet problem, which was given in the previous section for $\phi=0$ and which we denote by $(-\Delta)_0^{-1}$. Then (\ref{I.3a}) becomes
\[
u=(-\Delta)^{-1}_0 f(x,u,\nabla u)\quad.
\]
The idea behind this reformulation is to ultimately use a fixed point argument for solving the PDE. This is done in the same vein as for proving Cauchy-Lipschitz-Picard existence theorem
for ODEs of the form $\dot{y}=f(t,y(t))$. One writes the solution as an integral
\[
y(t)=y(t_0)+\int_{t_0}^tf(s,y(s))\ ds\quad,
\]
before applying a fixed point argument. The difficulty is to find the right space in which to work. In the ODE case, for smooth $f$, the integral operator
is very explicit and one shows that the space of Lipschitz functions gives an ad-hoc framework for the fixed point argument to work. In the case of PDEs, one must first have a thorough understanding of the operator $(-\Delta)_0^{-1}$ and how it acts on elements of various Banach spaces.\\

In two fundamental papers published in 1959 and in 1964, Schmuel Agmon, Avron Douglis, and Louis Nirenberg solved the elliptic Cauchy problems and the invertibility of elliptic operators of arbitrary orders in domains (see \cite{ADN1} and \cite{ADN2}). They worked in the context of Banach $L^p$ spaces. They obtained a series of optimal results that opened
the way to explore not only linear, but also non-linear PDEs, which were beforehand completely out of reach. We give below the simplest example of boundary problem solved by the Agmon-Douglis-Nirenberg 
theory of Cauchy problems in Banach spaces bearing in mind that the theory applies to a very wide family of problems.
 \begin{Th}
 \label{th-II.1}{\bf[Agmon, Douglis, Nirenberg 1959]} Let $\Omega$ be a smooth bounded domain of ${\R}^n$ and $p\in (1,+\infty)$. For any $f\in L^p(\Om)$, there exists
 a unique  solution $u$ to
 \be
\label{I.2}
\lf\{
\begin{array}{l}
\ds-\Delta u=f\quad\quad\mbox{ in }\Omega\\[5mm]
\ds u=0 \quad\quad\mbox{ on }\p\Omega\quad.
\end{array}
\rg.
\ee 
such that the distributional second and first derivatives of $u$ lie in $L^p(\Omega)$, and moreover
\[
\int_{\Omega}\sum_{|\al|\le 2}|\p^\al u|^p\le C(\Om,p)\ \int_\Om|f|^p\quad,
\]
where $C(\Om,p)$ is a positive constant independent of $f\in L^p(\Om)$.\hfill $\Box$
 \end{Th}
With such an estimate at hand, it becomes possible to use the operator $(-\Delta)_0^{-1}$ and a fixed point argument in the Banach spaces $L^p$ so as to obtain a solution to non-linear PDEs. A major piece of the puzzle was brought in in 1952 by Alberto Calder\'on and Antoni Zygmund. The``Singular integral theory'' says in particular\footnote{It is striking that the trace of the Hessian matrix encodes the whole Hessian matrix, in the $L^p$ space perspective.
} that, for {\it Schwartz functions}\footnote{The space of Schwartz function in ${\R}^n$ is the space of $C^\infty$ functions whose derivatives of arbitrary order decrease faster than any polynomial at infinity. It is preserved by the Fourier transform.}, $u$, the mapping
\[
\Delta u\longrightarrow (\p^2_{x_i x_j} u)_{i,j=1\cdots n}
\]
is continuous for the $L^p({\R}^n)$ norm into $L^p({\R}^n)$ for any $p\in (1,+\infty)$ (see \cite{CaZ1} and \cite{CaZ2}). The Agmon-Douglis-Nirenberg $L^p$-theory is parallel to a previous theory developed by the Polish mathematician Juliusz Schauder around 1930, which involves the more regular  {\it H\"older spaces}\footnote{These spaces were discovered by Otto Ludwig H\"older in his dissertation in 1882. He had been asked to characterize
the optimal regularity of the charge distribution $f$ that would ensure that the potential $u$ solving (\ref{I.2}) have continuous second derivatives. He discovered the existence of $\al>0$ such that
\be
\label{oo}
\sup_{x\ne y}\frac{|f(x)-f(y)|}{|x-y|^\al}<+\infty
\ee
is a sufficient condition implying that $u$ is $C^2$. This norm is denoted $\|f\|_{C^{0,\al}}$.} instead of {\it $L^p$ } or {\it Sobolev spaces}.

\section{Inequalities and a-priori estimates}
\reset
\subsection{Gagliardo-Nirenberg interpolation inequalities}
The field of PDEs is structured by {\it inequalities}. They are the ``power horses'' of the field. In the mid 1950s, functional analysis was already rich in inequalities: H\"older, Minkowski, Poincar\'e, Poincar\'e-Wirtinger, Young, Hausdorff-Young, Hardy, Hardy-Littlewood...etc and the more recent Sobolev inequalities. Deep scaling considerations (a common trait in Louis Nirenberg's work) led  Louis Nirenberg around 1959 (see \cite{Ni1} and \cite{Ni2}), and independently Emilio Gagliardo ( see \cite{Gag}), to discover a large family of inequalities for Sobolev norms ``sitting'' between two others. As an illustration, we have that
for any $p\in[1,+\infty]$, any $r\in [1,+\infty]$, and any $n\in {\N}^\ast$, there exists a constant $C>0$ such that for all Schwartz functions $u\in {\mathcal S}({\R}^n)$, there holds:
\be
\label{I.3}
\lf[\int_{{\R}^n}\sum_{|\al|= 1}|\p^\al u|^q\ dx^n\rg]^{1/q}\le\ C\ \lf[\int_{{\R}^n}|u|^p\ dx^n\rg]^{1/{2p}}\lf[\int_{{\R}^n}\sum_{|\al|\le 2}|\p^\al u|^r\ dx^n\rg]^{1/2r}
\ee
where $q^{-1}=2^{-1}\,(p^{-1}+r^{-1})$ and $dx^n$ denotes the Lebesgue measure on ${\R}^n$. In the family of {\it Gagliardo-Nirenberg interpolation inequalities}, one also finds the following one that we shall use later on. For any
$1\le p\le q<+\infty$ and any $n\in {\N}^\ast$, there exists $C>0$ such that for all Schwartz functions $u\in {\mathcal S}({\R}^n)$ there holds:
\be
\label{I.4}
\lf[\int_{{\R}^n}|u|^q\ dx^n\rg]^{1/{q}}\le\, C\, \lf[\int_{{\R}^n}|u|^p\ dx^n\rg]^{t/{p}}\ \lf[\int_{{\R}^n}\sum_{|\al|\le 1}|\p^\al u|^n\ dx^n\rg]^{(1-t)/n}
\ee
where $t=p\,q^{-1}$. The discovery of these estimates, together with the pioneering works of Marcel Riesz, Olaf Thorin, and J\'ozef Marcinkiewicz, lies at the origin of an important subfield of functional analysis developed by Jacques-Louis Lions, by Jaak Peetre (see \cite{LP}) and by Alberto Calder\'on (\cite{Cal}), and called {\it interpolation theory}.

\subsection{The use of Gagliardo-Nirenberg inequalities for proving a-priori estimates}

The notion of {\it a-priori} estimates is central in PDEs. Roughly speaking, it deals with finding an estimate of the form
\[
\|u\|_{E}<+\infty\quad,
\]
for some suitable Banach space $E$, all the while assuming we have a solution $u$ of some PDE, but prior to having actually proved the existence of such a solution. In concrete situations, looking at a given non-linear PDE problem, one establishes such an {\it a-priori} bound in order to perform one of the numerous available analytical methods in order to finally prove the existence of a solution satisfying that bound: fixed point argument in a perturbative approach, continuity method\footnote{This method, probably first used by Serge Bernstein in 1905,
consists in ``deforming'' continuously a given PDE problem to one which is directly solvable, in such a way as to preserve the {\it a-priori estimates}. The continuity method applied to elliptic non-linear problems is roughly done as follows. Once the {\it a-priori estimate} is established, one linearises the problem with respect to the deformation parameter , then calls upon the linear Agmon-Douglis-Nirenberg $L^p$ theory or the Schauder $C^{0,\al}$ theory, and finally uses the local inversion theorem to ensure the ``prolongation'' of the existence of the solution as the deformation parameter evolves.}, topological techniques (ex: Leray-Schauder theory), functional analysis approaches (Ex. monotone operator theory, Hille Yosida), successive approximation (ex : Galerkin method, convex integration), penalization approaches (ex : elliptization or viscosity method), variational appraoches (ex : minimization, min-max methods, Morse theory)... 

\medskip

Gagliardo-Nirenberg inequalities are mostly used to control non-linear terms in PDEs, and to establish {\it a-priori estimates}.
There are countless applications for these inequalities. In order to illustrate their might in dealing with non-linearities, the author
of these notes is cherry-picking a subject dear to his heart: the Dirichlet problem for maps taking values into submanifolds, also called {\it Harmonic map} or {\it vectorial Dirichlet problem}. This problem has many applications in geometry (in minimal surface theory and
in complex geometry -- it is used to describe the Teichm\"uller space of 2-dimensional surfaces, ...) as well as in physics (for instance the Dirichlet problem for maps taking values into the sphere is the
main mathematical object of the Ericksen-Leslie modeling of liquid crystals). 

\medskip\noindent
We consider exactly the same problem as the one posed by Dirichlet, that is to find critical points of the Dirichlet energy $E$ on a smooth bounded domain $\Omega$ for some boundary condition, but this time
under the additional constraint that the ``unknown'' $u$ take value in a given closed sub-manifold $N^n\subset {\R}^m$ of Euclidian space. Having fixed $\phi\in C^1(\p \Om, N^n)$, we seek a critical point $u$ of $E$ with the constraint that $u\in N^n$ 
everywhere. This constraint generates a new equation generalizing the Laplace equation (\ref{I.1}) and one shows that the problem is equivalent to
\be
\label{I.1b}
\lf\{
\begin{array}{l}
\ds\Delta u \in T_{u(x)}N^n\quad\quad\mbox{ in }\Omega\\[5mm]
\ds u=\phi \quad\quad\mbox{ on }\p\Omega\quad.
\end{array}
\rg.
\ee
where $T_{u(x)}N^n$ denotes the tangent space to $N^n$ in ${\R}^m$ at the point $u(x)$. The equation $\Delta u \in T_{u(x)}N^n$  is a natural generalization of the Laplace equation previously obtained for the unconstrained {\it Dirichlet problem}. Indeed, it says that the tangential component of $\Delta u$ to the constraint  is zero which simply implies that $\Delta u\equiv 0$ when there is no constraint. (\ref{I.1b}) can be recast as an elliptic non-linear equation of the form
\be
\label{I.2b}
\lf\{
\begin{array}{l}
\ds-\Delta u = f(u,\nabla u)\quad\quad\mbox{ in }\Omega\\[5mm]
\ds u=\phi \quad\quad\mbox{ on }\p\Omega\quad.
\end{array}
\rg.
\ee
This PDE is said to be {\it semi-linear}, because the term involving the highest order derivatives is linear. In contrast with the original Dirichlet problem (\ref{I.1}), one can show that (\ref{I.2b}) sometimes has more than one solution (and sometimes has infinitely many of them!). \\
To which extent is this non-linear Dirichlet problem degenerate, and what sufficient conditions would ensure that (\ref{I.2b}) is a well posed {\it Cauchy problem}?  Taking two solutions $u$ and $v$, the difference $u-v$ satisfies
\be
\label{I.3b}
u-v=(-\Delta)^{-1}_0\lf( f(u,\nabla u)-f(v,\nabla v)\rg)\quad.
\ee
In order to proceed to a ``contraction mapping argument'' we need to find a function space $E$ in which the difference of the non-linearities $(-\Delta)^{-1}_0\lf(f(u,\nabla u)-f(v,\nabla v)\rg)$ is controlled
by the norm of the difference $u-v$ in the same space, weighted by a constant $k<1$:
\be
\label{I.5b}
\|(-\Delta)^{-1}_0\lf(f(u,\nabla u)-f(v,\nabla v)\rg)\|_E\le k\ \|u-v\|_E\quad.
\ee
Such an inequality, also called {\it a-priori estimate} -- since we haven't yet settled the question of existence neither for $u$ nor for $v$ -- offers the possibility to ``absorb'' the non-linear right-hand side of (\ref{I.3b}) into the linear left-hand side $u-v$, and in turn deduce that $u=v$. \\

The {\it Gagliardo-Nirenberg interpolation inequalities} are essential to prove such {\it a-priori estimates} as (\ref{I.5b}). Our example, picked among countless others, was chosen to illustrate the possibility to solve Cauchy problems for non-linear PDEs. It is however representative of the might of the Gagliardo-Nirenberg inequalities.\\

Among the various inequalities pertaining to {\it vectorial Dirichlet problems}, there 
is a particularly elegant one\footnote{This inequality, which is a particular case of (\ref{I.4}), also appears in the work of Olga Ladyzenskaya.}, which, in combination with the Agmon-Douglis-Nirenberg $L^p$ theory, enables
to show the well-posedness of the {\it vectorial Dirichlet problem} for two dimensional domains under small energy assumptions. Namely, there exists a constant $C(\Om)>0$ such that for any smooth function $w$ supported in the two dimensional domain $\Omega$ the following inequality holds:
\be
\label{I.4b}
\int_{\Omega}|w|^4\ dx^2\le\, C(\Omega)\, \int_{\Omega}|w|^2\ dx^2\ \int_{\Omega}|\nabla w|^2\ dx^2\quad.
\ee
This inequality was centrally used by Michael Struwe in the framework of the {\it vectorial Dirichlet problem} to obtain {\it a-priori estimates} that imply the existence and uniqueness for the corresponding flow, also called {\it called harmonic map flow}, in two dimensions (see \cite{Str}). 

\section{The John-Nirenberg BMO space: when elasticity meets harmonic analysis}
\reset
The analysis of PDEs has evolved and keeps evolving in close ``partnership''  with the development of functional analysis and function space theory. Many linear and non-linear problems in PDE have stimulated the introduction of new function spaces -- we have already outlined above the importance
of the Sobolev spaces for solving the Dirichlet problem. The converse is also true:  knowledge and properties of certain function spaces can trigger a new understanding of PDE problems.\\

In 1961, the mathematician Fritz John was studying a rigidity problem from elasticity (\cite{Joh}). The strain exerted on a perfect elastic solid can be measured by the distance of the gradient of the resulting deformation with respect to the orthogonal group. In relation with this fundamental notion in elasticity, he asked the following question:
\[
\begin{array}{c}
  \mbox{\it Is it true that if the gradient of a transformation $f$ from Euclidian space ${\R}^n$}\\ 
  \mbox{\it  into itself  is ``close to'' the group of rotations at every point, } \\ 
  \mbox{\it  then it is globally close to one \underbar{single} rotation?} 
  \end{array}
\]
By ``close to'', it was originally meant in the $L^\infty$ norm. This rigidity question finds its origin in a work of the mathematician Arthur Korn from 1914 (\cite{Kor}). Extending Korn's results and now celebrated inequalities,
Kurt Friedrichs proved in an important work in 1947 (see \cite{Fri}) that the gradient of a deformation
is everywhere antisymmetric if and only if it is constant, provided that suitable boundary conditions,
which exclude rigid motions, are imposed. This result is an infinitesimal version of the question asked by John.

F. John gave first a relatively straightforward  counter-example to the $L^\infty$ version of the question, but he also proved that on every ball $B_r(x)\subset {\R}^n$ of arbitrary center $x$ and radius $r>0$, there exists a rotation $R_{x,r}$ such that
\[
\frac{1}{|B_r(x)|}\int_{B_r(x)}|\nabla f-R_{x,r}|\ dx^n\le C\ \|\mbox{dist}(\nabla f,SO(n))\|_{\infty}\quad,
\] 
where $|B_r(x)|$ is the volume of the ball $B_r(x)$. Thus John was proving that although the gradient of such a deformation could not be close to one single rotation globally, it is
 in average in the $L^1$-norm at any scale, close to a rotation $R_{x,r}$ that possibly depends on the ball. \\

In a subsequent collaboration \cite{JN}, which has since become a milestone in analysis, Fritz John and Louis Nirenberg systematically studied the sub-space of locally integrable functions,  called space of functions of
{\it Bounded Mean Oscillation} ($BMO$), whose elements satisfy
\[
\sup_{x\in{\R}^n,r>0}\frac{1}{|B_r(x)|}\int_{B_r(x)}|u-\ov{u}_{x,r}|\ dx^n<+\infty\quad,
\]
where $\ov{u}_{x,r}$ is the average of $u$ on the ball $B_r(x)$. They proved that this space is strictly larger than $L^\infty$, the space of globally bounded functions (for example, $\log|x|\in BMO\setminus L^\infty$), but that it is smaller than $L^p_{loc}$ for any $p<+\infty$. Precisely following an ingenious decomposition of ${\R}^n$ of {\it Calder\'on-Zygmund type} they proved the existence of $\al_n>0$ such that for any ball $B\subset{\R}^n$
\[
\int_{B} \exp\left[ \alpha_n\frac{|u-\overline{u}|}{\ \|u\|_{BMO}}\right]\ dx^n\le C_n\ |B|
\]
where $\ov{u}$ is the average of $u$ on $B$. The existence of such a bound for any ball is proved to characterize uniquely the space $BMO$.\\

The space $BMO$, which naturally arose in the context of elasticity in 1960, was apparently unknown to functional analysts. It was therefore a big surprise to discover, after the remarkable work of Elias Stein and Charles Fefferman in 1972 (\cite{FeSt}), that $BMO$ was the Banach dual of a famous space introduced in complex function theory some 40 years earlier by  Friedrich Riesz \cite{Rie}
and named ``Hardy space'' after a famous work by Godfrey Hardy from 1915 \cite{Har}. Historically, the Hardy space ${\mathcal H}^1$ was defined in the context of holomorphic functions on the disc $D^2$: it is made of the traces on the circle $S^1$ of holomorphic functions $f$ such that
\[
\lim_{r\rightarrow 1^-}\int_0^{2\pi} |f(r e^{i\,\theta})|\ d\theta<+\infty\quad.
\]
Note that this integral is an increasing function of the parameter $r$ for any holomorphic function $f$. The Hardy space was later extended to a much broader context of real function space theory. The dual spaces ${\mathcal H}^1$ and $BMO$ play a fundamental
role in PDEs. Empirically, one could say that they are the ``natural replacements'' for  $L^1$ and $L^\infty$. These two spaces are not compatible with Calder\'on-Zygmund theory, and in fact the Agmon-Douglis-Nirenberg results do not hold either for $L^1$ or for $L^\infty$. In contrast, ${\mathcal H}^1$ and $BMO$ are well-behaved in these theories.

It is unfortunately beyond the scope of this presentation to expose in its full glory the usefulness of the duality ${\mathcal H}^1-BMO$ in the analysis of PDEs. It plays a central role for instance in a theory called {\it integrability by compensation}, where some non-linear quantities appear in the form of products with such algebraic structures that improved integrability might be deduced. This has frequently been used in fluid mechanics as well as in geometry.\\

We content ourselves with mentioning one application. Using some straightforward integration by parts, one proves the following {\it Gagliardo Nirenberg interpolation inequality}: for any $n\in {\N}^\ast$ there exists $C_n>0$ such that for all Schwartz functions $u$ in ${\mathcal S}({\R}^n)$ one has
\[
\int_{{\R}^n}\sum_{|\al|= 1}|\p^\al u|^4\ dx^n\le\ C_n\ \|u\|^2_\infty\ \int_{{\R}^n}\sum_{|\al|= 2}|\p^\al u|^2\ dx^n\quad.
\]
Using more sophisticated arguments (Littlewood-Paley decomposition of tempered distributions) the author of these notes in collaboration with Yves Meyer prove a ``slight" improvement
of this inequality by replacing the $L^\infty$ norm of $u$ with its $BMO$-norm, namely: 
\[
\int_{{\R}^n}\sum_{|\al|= 1}|\p^\al u|^4\ dx^n\le\ C_n\ \|u\|^2_{BMO}\ \int_{{\R}^n}\sum_{|\al|= 2}|\p^\al u|^2\ dx^n\quad.
\]
This improved inequality entering in a larger class of inequalities by David Adams and Michael Fraizer (see \cite{AF}) has been crucially used in \cite{MR} for proving a partial regularity result for {\it stationary Yang-Mills fields}, which was obtained independently by Terence Tao and Gang Tian (see \cite{TT}).

 
\section{The maximum principle}
\reset
\subsection{Nirenberg's strong maximum principle for parabolic equations}
It would be impossible to speak about the work of Louis Nirenberg without mentioning the {\it maximum principle}. The contrast between the immense range of applications
of this principle and the simplicity of the heuristic idea behind it is amazing.\\

In one dimension, the maximum principle states that a continuously twice differentiable function on the segment $[0,1]$ satisfying
\[
u''\ge 0\quad\mbox{(resp. }u''\le 0)\quad\mbox{on }[0,1]
\] 
achieves its maximal (resp. minimal) value on the boundary of the segment, i.e. on $\{0\}\cup\{1\}$.

\noindent
In higher dimension, the maximum principle was known to Gauss since 1839 for solutions of the Laplace equation, owing to the mean value theorem for harmonic functions:
\[
\begin{array}{c}
\mbox{\it A solution to the Laplace equation on a bounded smooth domain }\\[1mm]
\mbox{\it achieves its extremal values on the boundary of the domain.}
\end{array}
\]  
This formulation requires $u$ to solve a PDE, and it is only at the beginning of the XXth century that, after successive contributions by Charles \'Emile Picard (1905) \cite{Pic}, Serge Bernstein (1910) \cite{Ber} and Leon Lichtenstein (1924) \cite{Lic}, the idea
of a general principle for elliptic partial differential \underbar{inequalities} emerged.\\

In a seminal five-page long paper published in 1927, \cite{Hop}, Eberhard Hopf opened the way to a wide range of applications of this principle by proving a general {\it strong version of the maximum principle}
for $C^2$ solutions to the following linear elliptic second order inequalities on a smooth bounded domain $\Om\subset {\R}^n$:
\be
\label{II.1}
L u:=\sum_{i,j=1}^n a_{ij}(x)\frac{\p^2u}{\p x_i\,\p x_j}+\sum_{i=1}^nb_i(x)\frac{\p u}{\p x_i}+c(x) u\ge 0\quad,
\ee
where $a_{ij}(x)$ is a map into uniformly positive definite symmetric matrices, and the coefficient $a_{ij}$ as well as $b_i$ and $c$ are bounded in $L^\infty(\Omega)$, with moreover $c(x)\le 0$ in $\Om$. A $C^2$ function $u$ satisfying (\ref{II.1}) is called a {\it sub-solution of $L$} (and a {\it super-solution of $L$} if $Lu\le 0$).\\

As in the 1-dimensional case, one formulates the {\it weak version of the maximum principle} as follows:
\[
Lu\ge 0\quad\mbox{in }\Om\quad\Longrightarrow\quad \max_{x\in \ov{\Om}}u(x)\le \max_{x\in \p{\Om}}u(x)\quad,
\]
whereas the {\it strong version of the maximum principle} discovered by Eberhard Hopf asserts that
\[
Lu\ge 0\quad\mbox{in }\Om\quad\mbox{and}\quad \exists \, x_0\in \Om\quad\mbox{s.t. }\max_{x\in \ov{\Om}}u(x)=u(x_0)\quad\Longrightarrow \quad u(x)\equiv u(x_0)\quad\mbox{ in }\Om\quad.
\]
In the early 1950s, a weak version of this principle was known to hold for {\it parabolic operators} of the type $\p_tu-Lu$. In 1953, Louis Nirenberg proved the corresponding strong version (see \cite{Ni3}).
\subsection{The notion of ``barriers'', the ``moving plane'' method, and the Gidas-Ni-Nirenberg symmetry principle}

The heuristic idea behind the strong maximum principle, at least in the simpler elliptic framework, has an interesting geometric representation. We say that a linear elliptic operator $L$ satisfies the strong maximum principle if the following holds. Let any pair of hyper-surfaces realized by two graphs of respectively a sub-solution $u$ and super-solution $v$ over a bounded
domain, with one of them sitting above the other (i.e. $u\ge v$) be given. The strong maximum principle says that if they touch at some \underbar{interior} point, then the two hyper-surfaces are necessarily identical. \\
A {\it sub-solution} satisfying $Lu\ge 0$ is then said to be a ``barrier'' with respect to a {\it super solution} satisfying $Lv\le 0$, and vice versa. This geometric interpretation is an incentive to manufacture {\it barriers} with respect to {\it solutions}, {\it sub-solutions}, and to {\it super-solutions} in order to prove {\it pointwise inequalities} via the maximum principle. 
This fruitful technique has become a classic in analysis, where it is known as a {\it comparison argument}. Devising suitable {\it barriers} is nothing short of being an art in itself. It requires deep intuition and thorough experience of the problems considered. \\

The geometric interpretation of the maximum principle in both its strong and weak formulations was probably first used in the 1955 work of Alexander Danilowitsch Alexandrov (see \cite{Ale} for a later publications of his original ideas that he presented for the first time in Z\"urich in 1955 according to Heinz Hopf \cite{Hopf}). He used a {\it comparison argument} between the solution itself and some reflections of it in order to prove that embedded constant mean curvature closed surfaces in ${\R}^3$ are necessarily isometric to a dilation of the unit sphere $S^2$. The relevance of the maximum principle and of elliptic theory in the resolution of this geometric problem is apparent in the equation satisfied by the graph $u$ of constant mean curvature $H$, namely: 
\[
\mbox{div}\lf[\frac{\nabla u}{\sqrt{1+|\nabla u|^2}}\rg]=2\, H\quad.
\]
Following an important paper by James Serrin \cite{Ser}, Louis Nirenberg in collaboration with Basilis Gidas and Wei Ming Ni ,  converted Alexandrov's original idea for constant mean curvature surfaces into a general method, as beautiful as it is efficient, nowadays known as the ``moving plane method'', see \cite{GNN}. With it, one can prove symmetry results (either with respect to a given direction, or full rotational symmetry) and uniqueness results for positive solutions to semi-linear scalar equations of the form
\be
\label{II.2}
-\Delta u=f(u)\quad.
\ee
These symmetry and uniqueness results are of utmost importance, since they extend to a non-linear framework a fundamental principle in quantum mechanics  and in spectral theory; stating that the ground state of the Laplace operator, which is necessarily positive,
\[
-\Delta u=\la_1\ u
\]
enjoys special symmetries and has multiplicity one (i.e. it is unique -- Krein-Rutman theorem).\\
This method and the ensuing symmetry results have important applications to diverse areas of science: the study of ground states of non-linear Schr\"odinger models in quantum mechanics, the vortex theory of Onsager in thermodynamics, turbulence in statistical physics, phase-transitions in Van der Walls fluids, the Yamabe problem in differential geometry (which is concerned with finding constant scalar curvature metrics in a given conformal class), etc. \\

The moving plane method consists in comparing an arbitrary solution $u$ for the semilinear equation (\ref{II.2}) with its successive reflections $u_\la$ across a continuous family of parallel hyper-planes. These reflections $u_\la$ are used as barrier functions for $u$. A key ingredient of the method is the strong version of the maximum principle and a refinement of it discovered
by E. Hopf in the 1950s, and now known as ``Hopf lemma''. It states that at maximum points on the boundary, the outward normal derivative of a {\it sub-solution} is strictly positive\footnote{This strict positivity of the normal outward derivative at the maximal point was already proved by Leon Lichtenstein in 1923 in \cite{Lic} for exact solutions of $Lu=0$ in 2 dimensions.} unless the {\it sub-solution} is identically constant.\\

Later on, Louis Nirenberg in collaboration with Henri Berestycki introduced a new method, still based on the {\it strong maximum principle}, and called ``sliding method'' (see for instance \cite{BN}). This new method was devised to prove various pointwise estimates and asymptotic behaviors for solutions in cylindrical domains to semi-linear equations of the form  (\ref{II.2}), as well as for their parabolic counterparts.
The {\it sliding method} has numerous important applications to traveling fronts problems in the mathematical modeling of combustion and flame propagation . Its novel idea consists in comparing the solution with its translations along the axis of the cylinder, rather than using the images by successive symmetries of the solution
as {\it barriers}.  It requires the use of a version of the maximum principle on ``narrow domains'' where the sign of the 0th order coefficient $c(x)$ is not required to be controlled due to S.R. Srinivasa Varadhan (see \cite{BNV}).

\subsection{The Dirichlet problem for non-linear second order elliptic equations}

We have stressed the importance of {\it a priori estimates} for solving the Dirichlet problem of semi-linear equations. In order to prove such estimates, we have combined the Agmon-Douglis-Nirenberg $L^p$-theory for boundary-value problems along with the Gagliardo-Nirenberg estimates in various Banach spaces. For many scalar equations of elliptic type which are more non-linear and also more degenerate than semi-linear equations, the maximum principle is an additional
tool that can be added into the mix to reach the desired estimates.\\

In a series of five fundamental papers written in collaboration with Luis Caffarelli and Joel Spruck, Louis Nirenberg identified the {\it maximum principle} as a fundamental device to obtain {\it a priori } estimates, and to solve the Dirichlet problem for  ``highly'' non-linear PDEs known as {\it fully non-linear PDEs} (see \cite{CNS1} \cite{CNS2} \cite{CNS3} \cite{CNS4} \cite{CNS5}). \\
An example of such an equation is the {\it Monge-Amp\`ere} equation
which appears in problems related to optimal transport as well as in geometric problems of prescribed curvatures. 
\begin{Th}
\label{th-II.2}
{\bf[Caffarelli, Nirenberg, Spruck 1984]}
Let $\Omega$ be a strictly convex domain. For any smooth and positive function $f$ in $\ov{\Omega}$, and for any smooth function $\phi$ in $\ov{\Om}$, there exists a unique strictly convex $C^\infty$ solution $u$
to
\[
\lf\{
\begin{array}{l}
\ds\mbox{det}\ \nabla^2 u=f\quad\quad\mbox{ in }\Om\\[2mm]
\ds u=\phi\quad\quad\mbox{ on }\p \Om\quad.
\end{array}
\rg.
\]
\end{Th}
The strict convexity requirement imposed on the solution forces ellipticity on the problem. To see this, let $F(\p^2 u):=\log\mbox{det}\,\nabla^2u$. Let $(u^{ij})_{i,j=1\cdots n}$ denote the inverse of the Hessian matrix of $u$. Assuming that $u$ is convex, the matrix $(u^{ij})$ is symmetric and positive definite. Hence: 
\[
\sum_{i,j}^n\frac{\p F}{\p r_{ij}}(r)\,\xi_i\,\xi_j=\sum_{i,j} u^{ij}\, \xi_i\,\xi_j>0\quad,
\]
which is tantamount to ellipticity. \\
The core of the argument involves several steps. First, one establishes {\it a-priori estimates} for the H\"older norm of the solution (assumed to exist). These estimates follow from the classical elliptic Agmon-Douglis-Nirenberg $L^p$-theory or Schauder theory for H\"older spaces, provided that $F(\p^2 u)$ is {\it uniformly elliptic}:
\[
\la\ |\xi|^2\le \frac{\p F}{\p r_{ij}}\ \xi_i\,\xi_j\le \Lambda\, |\xi|^2\quad,
\]
for some constants $0<\la<\La$, and also provided that the modulus of continuity of $\p_rF(\p^2 u)$ can be controlled up to the boundary. In other words, the task is to prove that $\p^2 u$ is uniformly bounded and uniformly continuous 
up to the boundary. This is done with the help of the maximum principle and comparison arguments. Eventually, the authors reach the {\it a-priori estimate}\footnote{This bound was independently obtained by Nikolai Vladimirovich Krylov \cite{Kry1}, using works of Lawrence Craig Evans \cite{Eva1}, \cite{Eva2} and Neil Trudinger \cite{Tru1}, \cite{Tru2}.}
\be
\label{ooi}
\|\p^2u\|_{C^{0,\al}(\ov{\Om})}\le C<+\infty\quad,
\ee
where $C^{0,\al}$ is the H\"older norm defined in (\ref{oo}). Once (\ref{ooi}) is established, the final part of the argument relies on a {\it continuity method}, where one interpolates $\phi$ with the determinant of the Hessian matrix of a
function which coincides with $\phi$ on the boundary. \\

In these highly intricate works, the {\it maximum principle} is showing its full quintessence and potential to obtain hidden {\it a-priori estimates}. The core of the argument involves the sophisticated construction of {\it barrier functions} with which the solution $u$ and its successive derivatives (or finite differences) are compared.\\

These five papers by Caffarelli, Nirenberg, Spruck and also by Joseph Kohn for one of them, have stimulated a tremendous amount of research activity on {\it fully non-linear equations} since their publications. These equations have an immense range of applications in many fields of science, including material sciences, finance, computer vision, ... The original ideas of Nirenberg {\it et al.} have influenced the development of a whole branch of analysis, called {\it viscosity theory for PDEs},  where the {\it maximum principle} plays a central and decisive r\^ole. The {viscosity theory for PDEs} was introduced by Lawrence Craig Evans \cite{Eva} and by Michael Crandall and Pierre-Louis Lions \cite{CL}.

\section{Solving problems from geometry}
\reset

The analysis of PDEs and differential geometry are intimately intertwined by essence. The central r\^oles played by the Laplace operator and by the $\ov{\p}$-operator in Riemann surface theory, constitutes the simplest illustration of this imbrication. The second half of the XXth century saw a dramatic acceleration of the transfer of techniques from non-linear PDEs to the resolution of problems that seem {\it a-priori} confined to geometry. A spectacular example of the might of the PDE approach in geometry is the recent proof of the Poincar\'e conjecture by Grigori Perelman, which heavily relies on the parabolic Ricci flow devised by Richard Hamilton.

\subsection{ Nirenberg's resolution of the Weyl problem}
The taste for geometry and the influence of geometric questions are manifest in Louis Nirenberg's work.  He is among the pioneers who introduced elaborate analysis tools for solving questions pertaining to embeddings, tensors, curvature, complex structures, ... His doctoral work itself dealt with geometry, and, following the invitation from his advisors James Stoker and Kurt Friedrichs, he solved a problem originally posed by Hermann Weyl:
\[
\begin{array}{c}
\mbox{\it Given a metric of positive Gauss curvature on the 2-dimensional sphere ,}\\[1.5mm]
\mbox{\it does there exist an isometric embedding of this sphere into a convex surface of }{\R}^3?
\end{array}
\]
Nirenberg answered this question assuming the given metric is four times continuously differentiable\footnote{Several mathematicians had already considered this question. Hans Lewy in the case of analytic metrics, and more generally Aleksei Pogorelov who solved the problem with the help of different methods.} in \cite{Ni0}. In this first work, the general ``philosophy'' we have discussed is already present, and it will remain recurrent in Louis Nirenberg's papers: one looks for {\it a-priori estimates}, and combine them with suitable {\it continuity methods} that leave the {\it a-priori estimates} unchanged along the deformation. The link with the previous section is made by writing the PDE satisfied by a graph $(x,u(x))$ of Gaussian curvature $K(x)$ at $(x,u(x))$: 
\[
\mbox{det}(\p^2 u)=K\ (1+|\p u|^2)^2\quad.
\]
Its associated Dirichlet problem was the subject of the last part of the first paper of the aforementioned series by L. Caffarelli, L. Nirenberg, and J. Spruck.

\subsection{The Nirenberg problem}

The original ``Nirenberg problem'' can be stated as follows:
\[
\begin{array}{c}
\mbox{\it For which function $K$ on $S^2$ does there exist a metric pointwise conformal}\\[2mm]
\mbox{ \it to the canonical metric $g_0$, and such that $K$ is the Gauss curvature of that metric?}
\end{array}
\]
By ``pointwise conformal'', we mean the existence of a function $u$ on $S^2$ such that $g=e^{2u}g_0$. Stated differently, the {\it Nirenberg problem} consists in establishing the existence (or lack thereof) of a solution to the Liouville equation
\[
-\Delta_{g_0}u=e^{2u}K-1\quad,
\]
where $\Delta_{g_0}$ is the {\it negative Laplace Beltrami} operator for $g_0$, the canonical metric of $S^2$. This simply formulated question has brought forth
an enormous amount of work since the early 1970s. Not only because it is the ``simplest" instance of a wide range of similar questions (higher dimension, different curvature tensor, scalar curvature, $Q$-curvature, $\sigma_k-$curvatures, fractional curvatures etc...), but also because it gives rise to the major issues faced by conformal geometric analysts in their study of ``critical non-linear PDEs'', such as {\it concentration of compactness phenomena}, {\it Palais-Smale sequences}, {\it Morse theory}, inter alia. These issues appear as well in many celebrated problems: the Yamabe problem, harmonic map theory, Yang-Mills equations, constant mean curvature surfaces, to name a few. The apparent simplicity of the Nirenberg problem fosters  the universal difficulties arising in conformal geometric analysis.\\

It would be much beyond the scope of the present report to give a detailed account of the various arguments and creativity which have flourished in the quest for solving Nirenberg's problem. We content ourselves with mentioning that not every choice for $K$ gives rise to a solution. This is seen, for example, using the
 Gauss-Bonnet theorem, and the now well-known Kazdan-Warner necessary condition. Let us also mention an important sufficient condition for the existence due to Alice Sun-Yung Chang and Paul Yang 
 \cite{ChY}. Let $K$ be positive on $S^2$ with only non-degenerate critical points and with at least two local maxima. Suppose further that $\Delta_{g_0} K>0$ at all saddle points of $K$. Then Nirenberg's problem has a solution.

\subsection{The Newlander-Nirenberg complex Fr\"obenius theorem}
Louis Nirenberg has made important contributions to complex geometry and complex analysis. Once again, PDE techniques lie at the heart of the approach he favored to tackle various geometric questions. \\
An important one deals with the {\it integrability of almost complex structures}. The question goes as follows. Let a map $J$ from $\R^{2n}$ into the space of real-valued $(n\times n)$ matrices be given. Assume it satisfies everywhere the condition
\[
J^2=-I_n\quad,
\] 
where $I_n$ is the $(n\times n)$ identity matrix. 
\[
\begin{array}{c}
\mbox{\it Does there exist a local diffeomorphism $w$ in a neighborhood of any point}\\[1.5mm]
\mbox{\it which transports the endomorphism $J$ into multiplication by $i$}\\[1.5mm]
\mbox{\it everywhere, after the canonical identification ${\R}^{2n}\simeq {\C}^n$?}
\end{array}
\]
For $n=1$, the question amounts to solving locally a Beltrami equation
\[
\p_{\ov{z}}w=\mu\ \p_{z} w\quad.
\]
This is successfully achieved by introducing the singular integral operator associated to the inverse of the Cauchy-Riemann operator $\ov{\p}^{-1}$ and by using an elementary fixed point argument. \\

 The general case $n>1$ is much more involved and leads to an overdetermined system of coupled linear Cauchy-Riemann type PDEs of the form
 \be
 \label{oao}
 \p_{\ov{z_j}}w=\sum_{k=1}^n\mu_{jk}\ \p_{{z_k}}w\quad,
 \ee
where $z_j$ are complex coordinates for the complex structure given by the value of $J$ at the point in the neighborhood of which we are working. The system being overdetermined (indeed, differentiating in $\ov{z_k}$ the $j$-th equation has to give the same result as differentiating in $\ov{z_j}$ the $k$-th equation, namely: $\p_{\ov{z_j}}\p_{\ov{z_k}}w=\p_{\ov{z_k}}\p_{\ov{z_j}}w$)
 there must be a structural constraint on the system for it to be {\it locally solvable}. A necessary condition was discovered by Albert Nijenhuis in his 1951 doctoral thesis, and independently by Paulette Libermann in \cite{Lib}, as well as by Beno Eckmann and Alfred Fr\"olicher in \cite{EF}\footnote{In these works, using the Cauchy-Kowalevski approach, this condition is also shown to be sufficient in the analytic framework.}. It can be stated as follows.
 We consider the complexification of the tangent space to ${\R}^{2n}$ at each point, and we call a $(0,1)$ vector any  vector of the complexified tangent space and of the form
 \[
 X-i\, J\,X\quad,
 \]
 where $X$ is a real vector. The necessary condition says that the space of $(0,1)$ vectors has to be stable under bracket operation. Louis Nirenberg and his student A. Newlander succeeded in proving that
 this {\it complex Fr\"obenius} type condition is in fact sufficient for the local solvability in the $C^\infty$ framework. \\
Note that, unlike in the one-dimensional case, applying the inverse of the Cauchy-Riemann operator $\p_{\ov{z_j}}$ to the right-hand side of (\ref{oao}) does not enable an iteration argument in any classical Banach space (H\"older, Sobolev). Indeed, an integration in $z_j$ is applied to a derivative in $z_k$ of $w$, for $k\ne j$: this leads nowhere. \\
Nirenberg and Newlander in \cite{NN} had the idea to recast the problem in terms of the inverse $z(w)$ of the diffeomorphism $w(z)$. It solves a PDE system of the form
 \[
 \p_{\overline{w_l}}z= F(w) \p_{{w}_l}z\quad.
 \]
Although this problem is now becoming non-linear (!), unlike the original one which was linear, it is this time possible to implement a scheme similar to the one devised in the one-dimensional case, since differentiation occurs with respect to only one of the independent variables in each equation. This problem is solvable in classical H\"older spaces.  

 \medskip
We have decided to end this presentation with this remarkable work, which can but trigger the admiration of any mathematician, even beyond the field of partial differential equation, and this over 58 years after its original publication!
 
 \newpage
\section{Conclusion}
\reset
At the end of these notes, one feels somewhat frustrated to have only presented one part of the prolific and monumental work of Louis Nirenberg. Many important contributions have been omitted such as the analyticity of solutions to analytic PDEs (in collaboration with Charles Morrey, \cite{MN}), the regularity of free-boundary problems (in collaboration with David Kinderlehrer and Joel Spruck, \cite{KiN1},\cite{KNS1}, \cite{KNS2}) or the partial regularity of solutions to the Navier-Stokes equation (in collaboration with Luis Caffarelli and Robert Kohn, \cite{CKN}) which to this day remains the optimal step towards solving the Millenium problem. 
The author of these notes also apologizes to
the numerous collaborators of Louis Nirenberg whose joint-works with him have not been mentioned here.\\
Louis Nirenberg's scientific endeavor is an exemplary reminder to all of us that research is first and foremost a collective venture, in which debating, discussing, and exchanging ideas play a decisive r\^ole. ..

It is the result of no coincidence that Louis Nirenberg made his professional home at the
Courant Institute at New York University. A prestigious institution that has fostered since its very creation a unique laboratory for the free exchange of scientific ideas. \\
Although there are still many important theoretical questions to answer\footnote{pertaining, inter alia, to very non-linear PDEs, degenerate PDEs, non-local PDEs, rough data, PDEs in connection with stochastics, and much more.}, the analysis of partial differential equations is nowadays mostly aimed at better understanding other fields of science, with applications in geometry, physics, mechanics, chemistry, biology, social sciences, technology. These developments, and the ones to come, anchor their roots on the immense efforts deployed in the last century by human intelligence in this area of mathematics.
Mathematical knowledge is however not only made of an accumulation of truths and results confined to papers and books, and transmitted in this form to future generations. A large and immaterial share of mathematical knowledge resides in the ``living part" of Mathematics, in mathematicians themselves, with their intuitions, their hesitations, their perseverance, and most importantly with their  esthetical quest and search for beauty (as surprising as it may sound to non-mathematicians!). Hermann Weyl once said
``{\it My work always tried to unite the truth with the beautiful, but when I had to choose one or the other, I usually chose the beautiful}''. \\
We do not know whether Louis Nirenberg would agree with this quote, but we would nonetheless like to thank him for the beautiful mathematics he has produced and for generously sharing it with us all for so many years.

\end{document}